\documentclass[a4paper,12pt]{amsart}
\usepackage{comment,amsmath,amssymb,amsthm,changepage,pifont,graphicx,tikz}
\usepackage[english]{babel}

\newcommand{\longhookrightarrow}{}
\DeclareRobustCommand{\longhookrightarrow}{\lhook\joinrel\relbar\joinrel\rightarrow}

\begin{document}

\title{A lower bound for the gonality conjecture}
\author{Wouter Castryck}
\date{}
\maketitle

\begin{abstract}
\vspace{-0.5cm}
For every integer $k \geq 3$ we construct a 
$k$-gonal curve $C$ along with a very ample divisor of degree $2g + k - 1$ (where $g$ is the genus of $C$) to which the vanishing statement from the Green-Lazarsfeld gonality conjecture does not apply.
\end{abstract}
\vspace{0.6cm}

The gonality conjecture due to Green and Lazarsfeld~\cite{greenlazarsfeld} states that for $C$ a smooth complex projective curve of genus $g \geq 1$ and gonality $k \geq 2$, and $L$ a globally generated divisor on $C$ of sufficiently large degree, one has the following vanishing criterion in Koszul cohomology:
\begin{equation} \label{gonalityconjecture}
 K_{i,1}(C,L) \neq 0 \quad \Leftrightarrow \quad 1 \leq i \leq h^0(C,L) - k - 1.
\end{equation}
This conjecture was proved two years ago by Ein and Lazarsfeld~\cite{einlazarsfeld}, who moreover provided the sufficient lower bound $\deg L \geq g^3$. In the meantime this was improved to $\deg L \geq 4g - 3$ by Rathmann~\cite{rathmann}.

It is likely that this bound can be improved further, although Green and Lazarsfeld already noted~\cite[p.\,86]{greenlazarsfeld} that one at least needs
\begin{equation} \label{effectivegonalityconjecture}
  \deg L \geq 2g + k - 1
\end{equation}
because of the non-vanishing of $K_{g-1,1}(K_C + D)$, for any divisor $D$ of rank $1$ and degree $k$ on $C$. Very recently Farkas and Kemeny~\cite{farkaskemeny} showed that if
$C$ is sufficiently generic inside the moduli space of $k$-gonal curves of genus $g$ then the bound \eqref{effectivegonalityconjecture} is sufficient for the vanishing criterion \eqref{gonalityconjecture} to hold. Moreover they conjectured that
this should be true for \emph{all} curves whose genus is large enough when compared to the gonality, a statement which they called the effective gonality conjecture. 
Results due to Green~\cite[Thm.\,3.c.1]{green} resp.\ Teixidor i Bigas~\cite[Prop.\,3.8]{teixidor} imply that this is indeed the case for trigonal resp.\ tetragonal curves of genus $g > 3$ resp.\ $g > 6$.

In this note we aim for an improved delimitation of the foregoing considerations, through the following result.\\

\noindent \textbf{Theorem.} \emph{For each $k \geq 3$ there exists a curve $C$ of genus $g = k(k - 1)/2$ along with a very ample divisor $L$ of degree $2g + k - 1$ such that
\begin{equation} \label{nonvanishing}
  K_{h^0(C,L) - k, 1}(C,L) \neq 0.
\end{equation}
In particular, the bound \eqref{effectivegonalityconjecture} is not sufficient for the gonality conjecture to apply.}\\

\noindent The main conclusions to be drawn are that on the one hand, at most, one can hope to improve Rathmann's bound to $\deg L \geq 2g + k$, and that on the other hand $g > k(k - 1)/2$ is a necessary lower bound in the statement of Farkas and Kemeny's  effective gonality conjecture. In particular one observes that the special cases implied by the works of Green and Teixidor i Bigas are sharp.

The construction is very short and explicit. Namely we let $C$ be a smooth projective plane curve of degree $k + 1$ and let $L$ be the effective divisor cut out by a curve of degree $k-1$. Note that
the gonality of $C$ equals $k$ by \cite[Prop.\,3.13]{serrano}
and that the degree of $L$ equals $(k-1)(k+1) = 2g + k - 1$, as announced.
This automatically entails that $L$ is very ample. Now the standard exact sequence
\[ 0 \longrightarrow \mathcal{O}_{\mathbf{P}^2}(-C) \longrightarrow \mathcal{O}_{\mathbf{P}^2} \longrightarrow \mathcal{O}_C \longrightarrow 0 \]
can be easily converted into
\[ 0 \longrightarrow \mathcal{O}_{\mathbf{P}^2}(-2) \longrightarrow \mathcal{O}_{\mathbf{P}^2}( k-1) \longrightarrow \mathcal{O}_C(L) \longrightarrow 0. \]
Taking cohomology and using that $H^0$ and $H^1$ vanish in the case of $\mathcal{O}_{\mathbf{P}^2}(-2)$ we end up with an isomorphism
\[ H^0(\mathbf{P}^2, \mathcal{O}_{\mathbf{P}^2}( k-1)) \longrightarrow H^0(C,L),  \]
which shows that $h^0(C,L) = {k+1 \choose 2}$ and that 
the embedding
\[ C \stackrel{|L|}{\longhookrightarrow} \mathbf{P}^{{k+1 \choose 2} - 1} \]
lands inside the $(k-1)$-th Veronese surface. 
 This implies that we have an injection
\[ K_{{k+1 \choose 2} - k, 1} \left(\mathbf{P}^2, \mathcal{O}_{\mathbf{P}^2}(k - 1) \right) \ \longhookrightarrow \ K_{{k+1 \choose 2}-k,1}(C,L) \ = \ K_{h^0(C,L) - k,1}(C,L).\]
But the former space is non-trivial, from which our theorem follows. One reference for this last fact is~\cite[App.]{green}, which can be applied to the decomposition $\mathcal{O}_{\mathbf{P}^2}( k-1) = \mathcal{O}_{\mathbf{P}^2}( k-2) \otimes \mathcal{O}_{\mathbf{P}^2}(1)$.
Alternatively we refer to~\cite[Thm.\,1.8]{bettisurfaces} for a more direct statement.\\
 

\noindent \emph{Remark.} We work over $\mathbf{C}$ because
for instance \cite{einlazarsfeld} does so as well, but the theorem presented above is valid over any algebraically closed field, with the same proof.

\section*{Acknowledgements}
The author was supported by the European Research Council
under the European Community's Seventh Framework Programme
(FP7/2007-2013) with ERC Grant Agreement 615722 MOTMELSUM,
by research project G093913N of the
Research Foundation Flanders (FWO), and by the Labex CEMPI (ANR-11-LABX-0007-01).
He wishes to thank an anonymous reviewer for simplifying the initial construction, and Gavril Farkas and Milena Hering for helpful suggestions and comments.

\small
\vspace{5mm}
\noindent \textsc{Laboratoire Paul Painlev\'e, Universit\'e de Lille-1}\\
\noindent \textsc{Cit\'e Scientifique, 59655 Villeneuve d'Ascq Cedex, France}\\
\vspace{-0.4cm}

\noindent \textsc{Departement Elektrotechniek, KU Leuven and imec-Cosic}\\
\noindent \textsc{Kasteelpark Arenberg 10/2452, 3001 Leuven, Belgium}\\
\vspace{-0.4cm}

\noindent \texttt{wouter.castryck@gmail.com}

\end{document}